\documentclass[12pt]{amsart}
\newtheorem{theorem}{Theorem}

\def\ds{\displaystyle}

\title[The Lempert function of the symmetrized polydisc]
{The Lempert function of the symmetrized polydisc in higher
dimensions is not a distance}

\author{Nikolai Nikolov, Peter Pflug and W\l odzimierz Zwonek}

\address
{Institute of Mathematics and Informatics\\ Bulgarian Academy of
Sciences\\ 1113 Sofia, Bulgaria}\email{nik@math.bas.bg}

\address{Carl von Ossietzky Universit\"at Oldenburg\\
Fachbereich Mathe\-ma\-tik\\ Postfach 2503\\ D-26111 Oldenburg,
Germany}\email{pflug@mathematik.uni-oldenburg.de}

\address{Instytut Matematyki, Uniwersytet Jagiello\'nski, Reymonta 4,
30-059 Krak\'ow, Poland}\email{Wlodzimierz.Zwonek@im.uj.edu.pl}



\begin{document}

\begin{thanks}
{This paper was written during the stays of the first and third
named authors at Universit\"at Oldenburg supported by grants from
the DFG (January -- March 2006 and November 2005 (DFG Projekt
227/8-1)). They like to thank both institutions for their support.
The third author was also supported by the Research Grant No. 1
PO3A 005 28, which is financed by public means in the programme
promoting science in Poland in the years 2005-2008.}
\end{thanks}

\begin{abstract} We prove that the Lempert function of
the symmetrized polydisc in dimension greater than two is not a
distance.
\end{abstract}

\maketitle

\section{Introduction}

A consequence of the fundamental Lempert theorem (see \cite{Lem})
is the fact that the Carath\'eodory distance and the Lempert
function coincide on any domain $D\subset\Bbb C^n$ with the
following property $(\ast)$ (cf. \cite{Jar-Pfl}):

\noindent ($\ast$) {\it  $D$ can be exhausted by domains
biholomorphic to convex domains.}

For more than 20 years it was an open question whether the
converse of the above result is true in the reasonable class of
domains (e.g. in the class of bounded pseudoconvex domains). In
other words, does the equality between the Carath\'eodory distance
and the Lempert function of a bounded pseudoconvex domain $D$
imply that $D$ has the property $(\ast)$.

The only counterexample so far, the so-called symmetrized bidisc
$\Bbb G_2,$ was recently discovered and discussed in a series of
papers (see \cite{Cos1}, \cite{Cos2}, \cite{Agl-You} and
\cite{Edi}, see also \cite{Jar-Pfl}).

What remained open is the following natural question (see
\cite{Jar-Pfl}):

{\it Do Carath\'eodory distance and Lempert function coincide on
the symmetrized polydisc $\Bbb G_n$ for any dimension $n\ge 3$?}

The aim of the present paper is to give a negative answer to the
above question proving that the Lempert function of $\Bbb G_n$
($n\ge 3)$ is not a distance. This implies that $\Bbb G_n$ ($n\ge
3)$ does not have property $(\ast)$ (for a direct proof of this
fact see \cite{Nik}).

Moreover, we show that for any dimension greater than two there
are bounded pseudoconvex domains not satisfying $(\ast)$ and for
which the Carath\'eodory distance and the Lempert function are
equal.

\section{Background and results}

Let $\Bbb D$ be the unit disc in $\Bbb C.$ Let
$\sigma_n=(\sigma_{n,1},\dots,\sigma_{n,n}):\Bbb C^n\to\Bbb C^n$
be defined as follows: $$\sigma_{n,k}(z_1,\dots,z_n)=\sum_{1\le
j_1<\dots<j_k\le n}z_{j_1}\dots z_{j_k},\quad 1\le k\le n.$$ The
domain $\Bbb G_n=\sigma_n(\Bbb D^n)$ is called {\it the
symmetrized $n$-disc}.

Recall now the definitions of the Carath\'eodory pseudodistance,
the Carath\'eodory-Reiffen  pseudometric, the Lempert function and
the Ko\-ba\-yashi-Royden pseudometric of a domain $D\subset\Bbb
C^n$ (cf. \cite{Jar-Pfl}): $$ \aligned
c_D(z,w)&=\sup\{\tanh^{-1}|f(w)|:f\in\mathcal O(D,\Bbb
D),f(z)=0\},\\ \gamma_D(z;X)&=\sup\{|f'(z)X|:f\in\mathcal O(D,\Bbb
D),f(z)=0\},\\ \tilde
k_D(z,w)&=\inf\{\tanh^{-1}|\alpha|:\exists\varphi\in\mathcal
O(\Bbb D,D): \varphi(0)=z,\varphi(\alpha)=w\},\\
\kappa_D(z;X)&=\inf\{\alpha\ge 0:\exists\varphi\in\mathcal O(\Bbb
D,D): \varphi(0)=z,\alpha\varphi'(0)=X\},
\endaligned
$$
where $z,w\in D,$ $X\in\Bbb C^n.$ The Kobayashi pseudodistance
$k_D$ (respectively, the Kobayashi--Buseman pseudometric
$\hat\kappa_D$) is the largest pseudodistance (respectively,
pseudonorm) which does not exceed $\tilde k_D$ (respectively,
$\kappa_D$).

It is well-know that $c_D\le k_D\le\tilde k_D,$
$\gamma_D\le\hat\kappa_D\le\kappa_D,$ and
$$\gamma_D(z;X)=\lim_{\Bbb C_\ast\ni t\to
0}\frac{c_D(z,z+tX)}{t}\ \hbox{(cf. \cite{Jar-Pfl})},$$ and if $D$
is taut, then
$$\ds\kappa_D(z;X)=\lim_{\Bbb C_\ast\ni t\to
0}\frac{\tilde k_D(z,z+tX)}{t}\ \hbox{(see \cite{Pang})}.$$

Repeat that for $m\in\Bbb N$,
$$k_D^{(m)}(z,w):=\inf\{\sum_{j=1}^m
\tilde k_D(z_{j-1},z_j): z=z_0, z_1,\dots,z_{m-1},z_m=w\in D\}.$$
Note that $\tilde k_D=k_D^{(1)}\ge k_D^{(2)}\ge\dots,$ $\ds
k_D=\lim_{m\to\infty}k_D^{(m)}$, and, if $D$ is taut, then
$$\hat\kappa_D(z;X)=\lim_{\Bbb C_\ast\ni t\to
0}\frac{k_D(z,z+tX)}{t}\ \hbox{(see \cite{KobM})}.\leqno{(1)}$$

For $m\in\Bbb N,$ consider the infinitesimal version of
$k_D^{(m)},$ namely
$$\kappa_D^{(m)}(z;X)=\inf\left\{\sum_{j=1}^m\kappa_D(z;X_j):
\sum_{j=1}^mX_j=X\right\}.$$ Then
$$\kappa_D=\kappa_D^{(1)}\ge\kappa_D^{(2)}\ge\dots\ge
\kappa_D^{(2n-1)}\ge\kappa_D^{(2n)}=\hat\kappa_D$$ (for last
equality see \cite{KobS}). We also point out that obvious
modifications in the proof of (1) in \cite{KobM} show that if $D$
is taut, then $$\lim_{u,v\to z,\ u\neq
v}\frac{k_D^{(m)}(u,v)-\kappa_D^{(m)}(z;u-v)}{||u-v||}=0$$
uniformly in $m$ and locally uniformly in $z;$ thus,
$$\ds\kappa_D^{(m)}(z;X)=\lim_{\Bbb C_\ast\ni t\to 0}\frac{
k_D^{(m)}(z,z+tX)}{t}$$ uniformly in $m$ and locally uniformly in
$z$ and $X.$

Note that $\Bbb G_n$ is a hyperconvex domain (see \cite{Edi-Zwo})
and, therefore, a taut domain. (Thus, all the introduced invariant
functions are continuous (in both variables) for $D=\Bbb G_n$.)
Even more, $\Bbb G_n$ is $c_{\Bbb G_n}$-finitely compact (see
Corollary 3.2 in \cite{Cos3}).

In the proof of our main result (Theorem 1) we shall need some
mappings defined on $\Bbb G_n.$

For $\lambda\in\overline{\Bbb D}$, $n\ge 2$ one may define the
rational mapping
$$
p_{n,\lambda}:\Bbb C^n\owns z=(z_1,\ldots,z_n)\mapsto(\tilde
z_1(\lambda),\ldots,\tilde z_{n-1}(\lambda))=\tilde z(\lambda)
\in\overline{\Bbb C}^{n-1},$$ where $\ds\tilde
z_j(\lambda)=\frac{(n-j)z_j+\lambda(j+1)z_{j+1}}{n+\lambda z_1},$
$1\le j\le n-1.$ Then $z\in\Bbb G_n$ if and only if $\tilde
z(\lambda)\in\Bbb G_{n-1}$ for any $\lambda\in\overline{\Bbb D}$
(see Corollary 3.4 in \cite{Cos3}).

We may also define for $\lambda_1,\ldots,\lambda_{n-1}\in\overline{\Bbb
D}$ the rational function
$$ f_{\lambda_1,\ldots,\lambda_{n-1}}=
p_{2,\lambda_1}\circ\ldots\circ p_{n,\lambda_{n-1}}:\Bbb
C^n\mapsto\overline{\Bbb C}.$$ Observe that
$$f_\lambda(z):=f_{\lambda,\dots,\lambda}(z)=\frac{\sum_{j=1}^{n}jz_j
\lambda^{j-1}}{n+\sum_{j=1}^{n-1}(n-j)z_j\lambda^j}.$$ By Theorem 3.2 in
\cite{Cos3}, $z\in\Bbb G_n$ if and only if
$\ds\sup_{\lambda\in\overline{\Bbb D}}|f_\lambda(z)|<1.$ In fact, by
Theorem 3.5 in \cite{Cos3}, if $z\in\Bbb G_n,$ then the last supremum is
equal to $\ds\sup_{\lambda_1,\dots,\lambda_{n-1}\in\overline{\Bbb D}
}|f_{\lambda_1,\dots,\lambda_{n-1}}(z)|$.

It follows that $$c_{\Bbb G_n}(z,w)\ge p_{\Bbb G_n}(z,w):=
\max_{\lambda_1,\dots,\lambda_{n-1}\in\Bbb T}|p_{\Bbb
D}(f_{\lambda_1,\dots,\lambda_{n-1}}(z),f_{\lambda_1,\dots,\lambda_{n-1}}(w))|,$$
where $\Bbb T=\partial\Bbb D$ and $p_{\Bbb D}$ is the Poincar\'e
distance; in particular,
$$\gamma_{\Bbb G_n}(0;X)\ge\lim_{\Bbb C_\ast\ni t\to
0}\frac{p_{\Bbb G_n}(0,tX)}{|t|} =\max_{\lambda\in\Bbb T}|\tilde
f_\lambda(X)|=:\rho_n(X),$$ where
$$\tilde f_{\lambda}(X)=\frac{\sum_{j=1}^{n}jX_j
\lambda^{j-1}}{n}.$$ Let $e_1,\dots,e_n$ be the standard basis of
$\Bbb C^n$ and $L_{k,l}=\{X\in\Bbb C^n:X=X_ke_k+X_le_l\},$ $1\le
k\le l\le n.$ Observe that if $X\in L_{k,l},$ then
$$\rho_n(X)=\frac{k|X_k|+l|X_l|}{n}.$$

For $n=2$ one has that $\kappa_{\Bbb G_2}\equiv c_{\Bbb G_2}\equiv
p_{\Bbb G_2}$ (see \cite{Agl-You,Cos1}). On the other hand, we
have the following.

\begin{theorem} Let $n\ge 3.$

(a) If $k$ divides $n,$ then $\ds\kappa_{\Bbb
G_n}(0;e_k)=\rho_n(e_k).$ Therefore, if $l$ also divides $n,$ then
$\kappa^{(2)}_{\Bbb G_n}(0;X)=\rho_n(X)$ for any $X\in L_{k,l}.$

(b) If $k$ does not divide $n,$ then $\ds\hat\kappa_{\Bbb
G_n}(0;e_k)>\rho_n(e_k).$

(c) If $X\in L_{1,n}\setminus(L_{1,1}\cup L_{n,n}),$ then
$\kappa_{\Bbb G_n}(0;X)>\rho_n(X).$

In particular, $k_{\Bbb G_n}(0,\cdot)\not\equiv p_{\Bbb
G_n}(0,\cdot),$ $\tilde k_{\Bbb G_n}(0,\cdot)\not\equiv
k^{(2)}_{\Bbb G_n}(0,\cdot),$ and $\Bbb G_n$ does not have
property $(\ast).$
\end{theorem}

\noindent{\bf Remarks.} (i) We already know that for $n\ge 3$ at
least one of the identities $\hat\kappa_{\Bbb
G_n}(0,\cdot)\equiv\gamma_{\Bbb G_n}(0,\cdot)$ and $\gamma_{\Bbb
G_n}(0,\cdot)\equiv\rho_n $ does not hold and, therefore, the same
applies to the identities $k_{\Bbb G_n}(0,\cdot)\equiv c_{\Bbb
G_n}(0,\cdot)$ and $c_{\Bbb G_n}(0,\cdot)\equiv p_{\Bbb
G_n}(0,\cdot).$ It will be interesting to know if however some of
them hold and whether $c^i_{\Bbb G_n}(0,\cdot)\equiv c_{\Bbb
G_n}(0,\cdot)$ ($c^i_{\Bbb G_n}$ denotes the inner Carath\'eodory
distance of $\Bbb G_n$).

(ii) Observe that ${\Bbb G_{2n}}_{|_{L_{n,2n}}}=\Bbb G_2.$ Then,
in contrast to (c), for $z,w\in L_{n,2n}$ one has that
$$p_{\Bbb G_{2n}}(z,w)\le\tilde k_{\Bbb G_{2n}}(z,w)\le \tilde k_{\Bbb
G_2}(z,w)=p_{\Bbb G_2}(z,w)\le p_{\Bbb G_{2n}}(z,w)$$ and
therefore $\tilde k_{\Bbb G_{2n}}(z,w)=p_{\Bbb G_{2n}}(z,w).$
\medskip

In spite of Theorem 1, for any $n\ge 3$ there are bounded
pseudoconvex domains $D\subset \Bbb C^n$ which do not have the
property $(\ast)$ but nevertheless $c_D\equiv\tilde k_D.$

\begin{theorem} Let $G\subset \Bbb C^m$ be a balanced domain
(that is, $\lambda z\in G$ for any $\lambda\in\overline{\Bbb D}$
and any $z\in G$). Then $D=\Bbb G_2\times G$ does not fulfill the
property $(\ast).$

On the other hand, if, in addition, $G$ is convex (for example,
$G$ is the unit polydisc or the unit ball), then $c_D\equiv\tilde
k_D.$

\end{theorem}

\section{Proofs}

\noindent{\it Proof of Theorem 1.} (a) We shall prove even more,
namely, that $\ds\kappa_{\Bbb G_n}(0;e_k)=\rho_n(e_k)$ if and only
if $k$ divides $n.$

Assume that $\kappa_{\Bbb G_n}(0;e_k)=\rho_n(e_k).$ Since $\Bbb
G_n$ is a taut domain, there exists an extremal mapping for
$\kappa_{\Bbb G_n}(0;e_k),$ that is, a holomorphic mapping
$\varphi:\Bbb D\to\Bbb G_n$ with
$\varphi(z)=(z\varphi_1(z),\dots,z\varphi_n(z)),$
$\ds\varphi_k(0)=1/\rho_n(e_k)=n/k$ and $\varphi_j(e_j)=0$ for
$1\le j\le n,$ $j\neq k.$ Observe that
$f_\lambda\circ\varphi\in\mathcal O(\Bbb D,\Bbb D)$ and
$f_\lambda\circ\varphi(0)=0$ for any $\lambda\in\overline{\Bbb
D}.$ It follows by the maximum principle that
$g_\lambda\in\mathcal O(\Bbb D,\overline{\Bbb D}),$ where
$$g_\lambda(z)=\frac{\sum_{j=1}^{n}j\varphi_j(z)
\lambda^{j-1}}{n+\sum_{j=1}^{n-1}(n-j)z\varphi_j(z)\lambda^j}.$$
Since $g_\lambda(0)=\lambda^{k-1},$ the maximum principle implies
that $g_\lambda\equiv\lambda^{k-1},$ that is
$$\sum_{j=1}^{n}j\varphi_j(z)\lambda^{j-1}=
n\lambda^{k-1}+z\sum_{j=1}^{n-1}(n-j)\varphi_j(z)\lambda^{k+j-1},\
\quad\lambda\in\overline{\Bbb D},z\in\Bbb D.$$ Comparing the
respective coefficients of these two polynomials of $\lambda,$ we
get that $\ds\varphi_k\equiv\frac{n}{k},$
$\varphi_1\equiv\dots\equiv\varphi_{k-1}\equiv
\varphi_{n+1-k}\equiv\dots\equiv\varphi_{n-1}\equiv 0$ and
$$(k+j)\varphi_{k+j}(z)\equiv(n-j)z\varphi_j(z),\ 1\le j\le n-k.$$
These relations imply that $\varphi_j\equiv 0$ if $k$ does not
divide $j,$ and $\ds\varphi_{j}\equiv\binom{n/k}{j/k}z^{j/k-1}$ if
$k$ divides $j.$ If $k$ does not divide $n,$ then $n-k<k[n/k]<n$
and hence $\varphi_{k[n/k]}\equiv 0,$ a contradiction. Conversely,
if $k$ divides $n,$ then put $\varphi=(\tilde
\varphi_1,\dots,\tilde\varphi_n),$ where $\tilde\varphi_j\equiv 0$
if $k$ does not divide $j$ and
$\ds\tilde\varphi_j(z)=\binom{n/k}{j/k}z^{j/k}$ if $k$ divides
$j.$ It follows from the proof above that $\varphi$ sends $\Bbb D$
into $\Bbb G_n,$ and, up to a rotation, it is the only extremal
mapping for $\kappa_{\Bbb G_n}(0;e_k)=\rho(e_k).$

To see that if $k$ and $l$ divide $n,$ then $\kappa^{(2)}_{\Bbb
G_n}(0;X)=\rho_n(X)$ for any $X\in L_{k,l},$ it is enough to
observe that
$$\rho_n(X)\le\kappa^{(2)}_{\Bbb
G_n}(0;X)\le\kappa_{\Bbb G_n}(0;X_ke_k)+\kappa_{\Bbb
G_n}(0;X_le_l)$$ $$=\rho_n(X_ke_k)+\rho_n(X_le_l)=\rho_n(X).$$

(b) Denote by $I,$ $J$, and $K$ the indicatrices of $\rho_n,$
$\hat\kappa_{\Bbb G_n}(0;\cdot)$, and $\kappa_{\Bbb
G_n}(0;\cdot),$ respectively ($I=\{X\in\Bbb C^n:\rho_{n}(X)<1\}$
and etc.). Note that if $X\in\overline{J}$ is an extreme point of
$\overline{I},$ then $X$ is an extreme point of $\overline{J}$ and
therefore $ X\in\overline{K}.$ Thus, (b) follows by the inequality
$\kappa_{\Bbb G_n}(0;ne_k/k)>1$  and the fact that $ne_k/k$ is an
extreme point of $\overline{I}.$ In fact, to see the last claim
observe that if $0<\alpha<1,$ $\rho_n(X)=\rho_n(Y)=1,$
$ne_k/k=\alpha X+(1-\alpha)Y$ and $\lambda\in\Bbb T,$ then
$$1=\lambda^{1-k}\tilde f_\lambda(ne_k/k)\le\alpha|\tilde
f_\lambda(X)|+(1-\alpha)|\tilde f_\lambda(Y)|$$
$$\le\alpha\rho_n(X)+(1-\alpha)\rho_n(Y)=1.$$ Hence, $\tilde
f_\lambda(X)=\tilde f_\lambda(Y)=\lambda^{k-1}$ for any
$\lambda\in\Bbb T,$ that is, $X=Y=ne_k/k.$

(c) First, note that if $\lambda\in\Bbb T,$ then the mapping
$(z_1,z_2,\dots,z_n) \to(\lambda
z_1,\lambda^2z_2,\dots,\lambda^nz_n)$ is an automorphism of $\Bbb
G_n$ and
$$\kappa_{\Bbb G_n}(0;\lambda X)=\kappa_{\Bbb
G_n}(0;X).$$ Applying these facts, we may assume that $X_1,X_n>0.$

Since
\begin{multline*}\kappa_{\Bbb G_n}(0;X)\ge\kappa_{\Bbb
G_{n-1}}(p_{n,1}(0);p_{n,1}'(0)(X))\\ =\kappa_{\Bbb
G_{n-1}}\left(0; \frac{n-1}{n}X_1e_1+X_ne_{n-1}\right),
\end{multline*}
 it follows by induction
on $n$ that $\kappa_{\Bbb G_n}(0;X)\ge\kappa_{\Bbb G_3}(0;Y),$
where $\ds Y=\frac{3X_1}{n}e_1 +X_ne_3.$ Assume that $\kappa_{\Bbb
G_n}(0;X)=\rho_n(X).$ Then $$\rho_n(X)\ge\ds\kappa_{\Bbb
G_3}(0;Y)\ge\rho_3(Y)=\rho_n(X)$$ and hence $\ds\kappa_{\Bbb
G_3}(0;Y)=\rho_3(Y).$ Now, taking an extremal mapping
$\varphi(z)=(z\varphi_1(z),z\varphi_2(z),z\varphi_n(z))$ for
$\kappa_{\Bbb G_3}(0;Y),$ with the same notations as in the proof
of (a), we obtain that $g_{\lambda}\in\mathcal O(\Bbb
D,\overline{\Bbb D}),$ $\lambda\in\Bbb T.$ Since $g_{\pm1}(0)=1,$
then $g_{\pm1}\equiv 1,$ that is
$$\varphi_1(z)\pm2\varphi_2(z)+3\varphi_3(z)=
3\pm2z\varphi_1(z)+z\varphi_2(z).$$ Thus,
$$\varphi_2(z)\equiv z\varphi_1(z)\hbox{ and
}\ds\varphi_3(z)\equiv 1+\frac{z^2-1}{3}\varphi_1(z).$$

Set $\psi=\varphi_1(z)/3.$ Then $g_\lambda\in\mathcal O(\Bbb
D,\overline{\Bbb D})$ means that
$$\left|\frac{\psi(z)+2\lambda z\psi(z)+\lambda^2(1+(z^2-1)\psi(z))}
{1+2\lambda z\psi(z)+\lambda^2z^2\psi(z)}\right|\le 1$$
$$\iff\left|\frac{\psi(z)(1+\lambda z)^2+\lambda^2(1-\psi(z))}{\psi(z)(1+\lambda z)^2
+1-\psi(z)}\right|\le 1$$
$$\iff\hbox{Re}(\psi(z)(1-\overline{\psi(z)})((\overline{\lambda}+z)^2-(1+\lambda
z)^2))\le 0.$$ If $\lambda=x+iy,z=iw,w\in\Bbb
R,a=\hbox{Re}(\psi(z))-|\psi(z)|^2,b=\hbox{Im}(\psi(z)),$ then
$$y(a(2w-y(w^2+1))+bx(1-w^2))\le 0,\ \forall\ x^2+y^2=1.$$ Setting
$x=0$ implies that $a\ge 0.$ Letting $y\to 0^+$ gives $-2aw\ge
(1-w^2)|b|.$ Hence $a=b=0$ if $w>0.$ Then the identity principle
implies that either $\psi\equiv0$ or $\psi\equiv1.$ Thus, either
$X_1=0$ or $X_n=0$ which a contradiction.\qed

\medskip

\noindent{\it Proof of Theorem 2.} The second part follows by the
equalities $c_{\Bbb G_2}=\tilde k_{\Bbb G_2}$ and $c_G=\tilde
k_G,$ and the product property of $c_D$ and $\tilde k_D:$
$c_D=\max\{c_{\Bbb G_2},c_G\}$ and $\tilde k_D=\max\{\tilde
k_{\Bbb G_2},\tilde k_G\}$ (cf. \cite{Jar-Pfl}).

The proof of the first part the proof in \cite{Edi} that $\Bbb
G_2$ does not have property $(\ast).$ For convenience of the
reader, we include it.

Let $h_1(z_1+z_2,z_1z_2)=\max\{|z_1|,|z_2|\},$ $\ds
h_2(z)=\inf\{t>0:z/t\in G\}$ (the Minkowski function of $G$),
$h=\max\{h_1,h_2\}$ and
$$\pi_{\lambda}(z_1,\dots,z_{m+2})=(\lambda z_1,\lambda^2 z_2,\lambda
z_3,\dots,\lambda z_{m+2}),\ \lambda\in\Bbb C.$$ Note that
$h(\pi_{\lambda}(z)=|\lambda|h(z)$ and $D=\{z\in\Bbb
C^{m+2}:h(z)<1\}.$

Assume now that $D$ fulfills the property $(\ast).$ Take two
points $a,b\in \Bbb G_2\times\{0\}\subset D.$ We may find an
$\varepsilon>0$ and a domain
$D_\varepsilon\subset\{h<1-\varepsilon\}$ which is biholomorphic
to a convex domain $\tilde D_\varepsilon$ and such that $\lambda
a,\lambda b\in D_\varepsilon$ for $\lambda\in\overline{\Bbb D}.$
Let $\varphi_\varepsilon:D_\varepsilon\to\tilde D_\varepsilon$ be
the corresponding biholomorphic mapping. We may assume that
$\varphi_\varepsilon(0)=0$ and
$\varphi'_\varepsilon(0)=\hbox{id}.$ Note that
$$g_\varepsilon(\lambda)=\varphi_\varepsilon^{-1}
\left(\frac{\varphi_\varepsilon(\pi_\lambda(a))+
\varphi_\varepsilon(\pi_\lambda(b))}{2}\right),$$ is a holomorphic
mapping from a neighborhood of $\overline{\Bbb D}$ into $D.$ We
have $$g_\varepsilon(0)=0,
g_{\varepsilon,1}'(0)=\frac{a_1+b_1}{2},g_{\varepsilon,_2}'(0)=0,\dots,
g_{\varepsilon,m+2}'(0)=0$$
$$\hbox{and }g_{\varepsilon,2}'(0)=a_2+b_2+\frac{c_\varepsilon}{4}(a_1-b_1)^2,
\hbox{ where
}c_\varepsilon=\frac{\partial^2\varphi_{\varepsilon,2}}{\partial
z_1^2}(0).$$ Thus, the mapping
$f_\varepsilon(\lambda)=\pi_{1/\lambda}\circ
g_\varepsilon(\lambda)$ can be extended at $0$ as
$$f_\varepsilon(0)=\left(\frac{a_1+b_1}{2},\frac{a_2+b_2}{2}+\frac{c_\varepsilon}
{8}(a_1-b_1)^2,0,\dots,0\right).$$ Since $D$ is assumed to satisfy
the property $(\ast),$ it should be pseudoconvex, that is, $h$ is
a plurisubharmonic function. Then the maximum principle implies
that $\ds
h(f_\varepsilon(0))\le\max_{|\lambda|=1}h(f_\varepsilon(\lambda))<1$
which means that $f_\varepsilon(0)\in D.$

Assuming now that $\ds\lim_{\varepsilon\to 0^+}c_\varepsilon\neq
0$ and having in mind that $c_\varepsilon$ is bounded, we may find
$c\neq 0$ such that
$$m=\left(\frac{a_1+b_1}{2},\frac{a_2+b_2}{2}+c(a_1-b_1)^2,0,
\dots,0\right)\in\overline{D}$$ for any $a,b\in\overline{\Bbb
G_2}\times\{0\}.$ Taking $\alpha=e^{i(\arg(c)+\pi)/2},$
$a_1=\alpha+1,$ $a_2=\alpha,$ $b_1=\alpha-1,$ $b_2=-\alpha,$ we
obtain that $\ds 1\ge h_1(m)=\frac{1+\sqrt{1+16c}}{2}$ which is
impossible.

Thus, $\ds\lim_{\varepsilon\to
0^+}f_\varepsilon(0)=\frac{a+b}{2}\in\overline{D}$ for any $a,b\in
\Bbb G_2\times\{0\},$ that is, $\Bbb G_2$ is a convex domain, a
contradiction (for example, $(2,1),(2i,-1)\in\partial\Bbb G_2$,
but $(1+i,0)\not\in\overline{\Bbb G_2}$).

\qed

\end{document}